\documentclass[a4paper, oneside, 11pt]{article}
\usepackage[english]{babel}
\usepackage{amsfonts}
\usepackage{amssymb}
\usepackage{graphics}
\usepackage{amsrefs}
\usepackage{latexsym}
\begin{document}

\title{{\bf{\Large{Bayesian nonparametric analysis for a species sampling model with finitely many types}}}\footnote{{\it AMS (2000) subject classification}. Primary: 60G58. Secondary: 60G09.}}
\author{\textsc {Annalisa Cerquetti}\footnote{Corresponding author: Annalisa Cerquetti, Dept. Studi Geoeconomici, Facolt\`a di Economia,  Sapienza Universit\`a di Roma, Via del Castro Laurenziano, 9, 00161, Rome, Italy. E-mail: {\tt annalisa.cerquetti@uniroma1.it}}\\
  \it{\small Department of Geoeconomical, Linguistic, Statistical and Historical Studies},\\ \it {\small Sapienza University of Rome, Italy }}
\newtheorem{teo}{Theorem}
\date{}
\maketitle{}

\begin{abstract}
We derive explicit Bayesian nonparametric analysis for a species sampling model with finitely many types of Gibbs form of type $\alpha= -1$ recently introduced in Gnedin (2009). Our results complement existing analysis under Gibbs priors of type $\alpha \in [0, 1)$ proposed in Lijoi et al. (2008). Calculations rely on a groups sequential construction of Gibbs partitions introduced in Cerquetti (2008).  
\end{abstract}

\section{Introduction}

In the {\it species sampling} problem a random sample is drawn from an hypothetically infinite population of individuals to make inference on the unknown total number of different species. A Bayesian nonparametric approach to this problem has been recently proposed in Lijoi et al. (2007, 2008) to derive posterior predictive inference on species richness for an additional sample under the assumption that the vector of multiplicities of different species observed is a random sample from an {\it exchangeable Gibbs} partition of type $\alpha \in [0,1)$. 

Exchangeable Gibbs partitions of type $\alpha \in [-\infty, 1)$ (Gnedin and Pitman, 2006) are models for random partitions of the positive integers which extend the Ewens-Pitman two-parameter family and are characterized by a probability function (EPPF) with the following structure
\begin{equation}
\label{gibbsGP}
p(n_1,\dots, n_k)=V_{n,k} \prod_{j=1}^k (1-\alpha)_{n_j-1 \uparrow},
\end{equation}
with weights $(V_{n,k})$ being solution to the {\it backward} recursion  $V_{n,k}=(n-\alpha k)V_{n+1, k}+V_{n+1,k+1}$  with $V_{1,1}=1$. Gnedin and Pitman even provide a constructive result for this class (cfr. Th. 12) showing that each element arises as a unique probability mixture of extreme partitions which are\\\\
{\it \begin{tabular}{ll}
a) $PD(\alpha, \xi|\alpha|)$ partitions with $\xi=1,\dots,\infty$ & for $\alpha \in [-\infty, 0)$,\\
b) $PD(0,\theta)$ partitions with $\theta \in [0,\infty)$ &  for $\alpha=0$,\\
c) $PK(\rho_\alpha|t)$ partitions with $t \in [0,\infty)$ & for $\alpha \in (0,1)$.
\end{tabular}\\
}\\\\
Here $PD(\cdot, \cdot)$ stands for the two-parameter Poisson-Dirichlet distribution (Pitman and Yor, 1997) and $PK(\rho_\alpha|t)$ for the conditional Poisson-Kingman distribution derived from the stable subordinator, (cfr. Pitman, 2003).\\\\
Lijoi et al. (2008) establish general distributional results and properties for EPPFs belonging to subclass c) to make conditional predictions according to a Bayesian nonparametric procedure. Focusing on class c) their analysis relies on the hypothesis that the total unknown number of species is so large to be assumed infinite. An assumption that may be unrealistic in concrete applications. Here we focus on a class of random partitions with finite but random number of different species, recently introduced in Gnedin (2009), which belongs to subclass a) and contribute, in view of possible future Bayesian implementations, deriving posterior predictive results analogous to that in Lijoi et al. (2008).\\

First recall that for $\alpha < 0$, $\theta=|\alpha|\xi$ and $\xi=1,2,3,\dots$  $PD(\alpha, \xi|\alpha|)$ model has EPPF (see e.g. Pitman, 1996, 2006)
\begin{equation}
\label{finitegibbs}
p(n_1, \dots, n_k)=\frac{(\xi|\alpha|- |\alpha|)_{k-1 \downarrow |\alpha|}}{(\xi|\alpha| +1)_{n-1 \uparrow 1}} \prod_{j=1}^{k} (1 +|\alpha|)_{n_j -1\uparrow}
\end{equation}
or equivalently
$$
p(n_1, \dots, n_k)=\frac{(\xi)_{k \downarrow}}{(|\alpha|\xi)_{n \uparrow 1}}\prod_{j=1}^{k}(|\alpha|)_{n_j \uparrow},
$$
and arises by the following sequential procedure. Given the partition of $[n]$ in $K_n= k$ blocks with occupancy counts ${\bf n}= (n_1, \dots, n_k)$, the partition of  $[n+1]$ is obtained  through one-step prediction rules 
$$
p_j({\bf{n}})=\frac{n_j  + |\alpha|}{n + |\alpha|\xi}    \mbox {\hspace{1.5cm}and {\hspace{1.5cm}}}      
p_{0}({\bf{n}})=\frac{|\alpha|(\xi-k)}{\xi|\alpha|+n}
$$
where $p_j({\bf n})$ for $j=1,\dots,k$ stands for the probability of randomly placing $n+1$ in an {\it old} block $j$, while $p_0({\bf n})$ stands for the probability of $n+1$ to form a new block $k+1$. 
For $\alpha = -1$ the EPPF in (\ref{finitegibbs}) reduces to   
$$
p(n_1, \dots, n_k)=\frac{(\xi-1)_{k-1\uparrow -1}}{(\xi+1)_{n-1}}\prod_{j=1}^{k} {n_{j\uparrow}}
$$
or alternatively
$$
p(n_1,\dots, n_k)=\frac{(\xi)_{k \downarrow}}{(\xi)_{n \uparrow 1}}\prod_{j=1}^{k}{n_{j\uparrow}}
$$
with one-step prediction rules
$$
p_j({\bf n})=\frac{n_j  + 1}{n + \xi}
\mbox{\hspace{1.5cm}   and  \hspace{1.5cm}     }
p_0({\bf n})=\frac{( \xi- k)}{n+\xi}.
$$ 

The corresponding limit frequencies $(\tilde{P}_{\xi, 1},\dots, \tilde{P}_{\xi, \xi})$ have a stick-breaking representation
$$
\tilde{P}_{\xi,j}=W_j\prod_{i=1}^{j-1}(1 -W_i), \mbox{    with independent     } W_i \sim Beta(2, \xi -i) 
$$
for $i=1\dots, \xi$ and Beta (2,0) a Dirac mass at $1$.

\section{Gnedin's model with finitely many types}

As from Gnedin and Pitman result, $PD(\alpha, \xi|\alpha|)$ models are extreme points of a convex set of Gibbs partitions of type $\alpha <0$, whose elements are in one to one correspondence with a set of mixing distributions over the set of the positive integers.  Gnedin (2009) studies the particular model arising  by mixing a $PD(\alpha, |\alpha|\xi)$ model for $\alpha =-1$ over $\xi$ with 
$$
P(K=\xi)=\frac{\gamma(1-\gamma)_{\xi-1}}{\xi!}
$$
for $\xi=1,2,\dots$ and $\gamma \in (0,1)$ and shows it has weights 
\begin{equation}
\label{pesi}
V_{n,k}=\frac{(k-1)!}{(n-1)!}\frac{(1-\gamma)_{k-1} (\gamma)_{n-k}}{(1+\gamma)_{n-1}}=\frac{(k-1)!}{(n-1)!}\frac{\gamma(1-\gamma)_{k-1} }{(\gamma+n-k)_{k}},
\end{equation}
hence EPPF
\begin{equation}
\label{gnedin}
p(n_1, \dots, n_k)= \frac{(k-1)!}{(n-1)!}\frac{(1-\gamma)_{k-1} (\gamma)_{n-k}}{(1+\gamma)_{n-1}} \prod_{j-1}^k n_j!
\end{equation}
obtained by sequential construction with one-step prediction rules
$$
p_{j}({\bf n}) =\frac{(n-k+\gamma)(n_j +1)}{n(n +\gamma)}          \mbox{ \hspace{0.5cm} for $j=1,\dots, k$  \hspace{1cm} and  \hspace{0.5cm} } p_0({\bf n})=\frac{k(k-\gamma)}{n(n+\gamma)}.
$$

Gnedin also derives analogous of the Ewens sampling formula and further results on the vector of frequencies and of exchangeable sequences induced by sampling from this model. \\\\

For what follows it is worth to notice that Gnedin's one-step prediction rules may be equivalently expressed as $m$-steps prediction rules as in Cerquetti (2008, Prop. 3), which here we formulate as a group sequential random allocation of balls labelled $1,2...$ in a series of boxes.  
First from (\ref{pesi}) we obtain
$$
V_{n+m,k+k^*}=\frac{(k+k^*-1)!}{(n+m-1)!}\frac{(1-\gamma)_{k+k^*-1} (\gamma)_{n+m-k+k^*}}{(1+\gamma)_{n+m-1}},
$$
then by basic properties of rising factorials specialization of the general formulas for Gibbs partitions easily follow. Start with box $B_{1,1}$ with a single ball. Given the placement of the first group of $n$ balls in a $(n_1, \dots, n_k)$ configuration in $k$ boxes, the {\it new} group of $m$ balls labelled \{n+1,\dots, n+m\} is: \\\\
a) allocated in the {\it old} $k$ boxes in configuration $(m_1,\dots, m_k)$, for $m_j \geq 0$, $\sum_{j=1}^k m_j=m$, with probability 
\begin{equation}
\label{alloldA}
p_{\bf m}({\bf n})=
\frac{1}{(n)_m}\frac{(\gamma +n -k)_m}{(\gamma +n)_m} \prod_{j=1}^k (n_j +1)_{m_j}
\end{equation}\\
b) allocated in $k^*$ {\it new} boxes in configuration $(s_1, \dots, s_{k^*})$, for $\sum_{j=1}^{k^*} s_j =m$, $1 \leq k^* \leq m$, $s_j \geq 1$, with probability  
\begin{equation}
\label{allnewA}
p_{\bf s}({\bf n})=\frac{(k)_{k^*}}{(n)_m}\frac{(k-\gamma)_{k^*} (\gamma+n-k)_{m-k^*}}{(\gamma +n)_{m}} \prod_{j=1}^{k^*} s_j!
\end{equation}\\
c) $s < m$ balls are allocatedat $k^*$ {\it new} boxes in configuration $(s_1,\dots,s_{k^*})$ and the remaining $m-s$ balls in the {\it old} boxes in configuration $(m_1,\dots, m_k)$ for $\sum_{j=1}^{m} m_j= m-s$, $1 \leq s \leq m$, $\sum_{j=1}^{k^*} s_j=s$, $m_j \geq 0$, $s_j \geq 1$ with probability
\begin{equation}
\label{oldandnew}
p_{{\bf s, m}}({\bf n})=\frac{(k)_{k^*}}{(n)_m}\frac{(k-\gamma)_{k^*} (\gamma+n-k)_{m-k^*}}{(\gamma +n)_{m}} \prod_{j=1}^{k}(n_j+1)_{m_j}\prod_{j=1}^{k^*} s_j!
\end{equation}

These $m$-steps prediction rules allows to readily obtain Bayesian posterior predictive distributional results for the random partition induced by an additional $m$-sample from Gnedin's model. By exploiting the definition of central and non-central Lah numbers as particular case of central and non-central generalized Stirling numbers of the first kind, the results of next section are obtained specializing results in Lijoi et al. (2008) by means of expressions derived in Cerquetti (2008). See the Appendix for the relationship between generalized Stirling numbers and generalized factorial coefficients both central and non-central.  For an explicit example of application of this kind of results see e.g. Section 4. in Lijoi et al. (2008).\\

\section{Posterior predictive analysis of Gnedin's model}

First notice that generalized Stirling numbers of the first kind $S_{n,k}^{-1, -\alpha}$ for $\alpha =-1$ admit an explicit expression known as {\it Lah numbers} 
$$
S_{n,k}^{-1, 1}= {n -1 \choose k-1} \frac{n!}{k!}
$$
which are connection coefficients defined by
$$
(x)_{n \uparrow} =\sum_{k=0}^n \frac{n!}{k!} {n  -1 \choose n-k}(x )_{k \uparrow -1}
$$
hence, by an application of (10) in Gnedin and Pitman (2006), (see Gnedin, 2009, eq. (2)) the distribution of the number $K_n$ of occupied boxes for Gnedin's model is readily calculated as
$$
Pr(K_n=k)=V_{n,k} S_{n,k}^{-1,1} = {n \choose k} \frac{(1-\gamma)_{k-1} (\gamma)_{n-k}}{(1+\gamma)_{n-1}}.
$$

Then recall the definition of non-central Lah numbers with parameter of non-centrality $r$ (see e.g. Charalambides, 2005) which correspond to non central generalized Stirling numbers for $\alpha = -1$
$$
S_{n,k}^{-1,1, r} =\frac {n!}{k!} {n -r -1\choose n-k}.
$$

Now a Bayesian nonparametric posterior predictive analysis of the class (\ref{gnedin}) is readily derived. Notice that for the sake of generality here we mantain the treatment in terms of random allocation of balls in boxes. The obvious traslation in terms of random partition of individuals among species follows easily.\\

Given the sufficiency of the number $K_n$ of boxes induced by the basic sample, the joint distribution of the number $S$ of balls allocated in new $k^*$ boxes in a specific configuration $(s_1, \dots, s_{k^*})$ given $(n_1, \dots, n_k)$ is obtained marginalizing (\ref{oldandnew}) with respect to $(m_1,\dots, m_k)$ and by an application of the multinomial theorem for rising factorials (see the Appendix) as in Lijoi et al. (2008) (cfr. eq. (27) in Cerquetti, 2008),
\begin{equation}
\label{s1sk}
Pr(s_1,\dots, s_{k^*}|K_n=k)=\frac{(k)_{k^*}}{(\gamma +n)_{m}}\frac{(k -\gamma)_{k^*} (\gamma +n -k)_{m-k^*}}{(n)_m}  {m \choose s}(n+k)_{m-s\uparrow}\prod_{j=1}^{k^*}s_j!.
\end{equation}

The joint distribution of the number of new boxes $K_m$ and of the total number of balls falling in new boxes $S$ given $K_n$ (cfr. eq. (28) in Cerquetti, 2008) is given by
\begin{equation}
\label{jointks}
Pr(K_m=k^*, S=s|K_n=k)=\frac{(k)_{k^*}}{(n)_m}\frac{(k -\gamma)_{k^*} (\gamma +n -k)_{m-k^*}}{(\gamma +n)_{m}}  {m \choose s}(n+k)_{m-s\uparrow} {s \choose k^*}\frac{(s-1)!}{(k^*-1)!}
\end{equation}
and arises by (\ref{s1sk}) by summing over the space of all partitions of $s$ elements in $k^*$ blocks and exploiting the definition of Lah numbers.\\

Marginalizing (\ref{jointks}) with respect to $K_m$, a probability distribution for the total number $S$ of balls in new boxes is obtained in terms of Lah numbers according to eq. (29) in Cerquetti (2008) and eq. (11) in Lijoi et al. (2008)
\begin{equation}
\label{laws}
Pr(S=s|n_1, \dots, n_k)={m \choose s}\frac{(n+k)_{m-s\uparrow}}{(n)_m(\gamma +n)_m}\sum_{k^*=0}^{s} {s \choose k^*} \frac{(s-1)!}{(k^*-1)!} (k)_{k^*} (k -\gamma)_{k^*} (\gamma +n -k)_{m-k^*}.
\end{equation}

The probability distribution of the number $K_m$ of new blocks induced by the additional sample given the basic sample follows marginalizing (\ref{jointks}) with respect to $S$ and exploiting the definition of non-central Lah numbers with parameter of non centrality $r=-(n+k)$. An application of eq. (4) in Lijoi et al. (2007) yields
\begin{equation}
\label{bioalpha}
Pr(K_m=k^* | K_n=k)=\frac{(k-\gamma)_{k^*}(\gamma +n -k)_{m-k^*}}{(n+\gamma)_m}{m \choose k^*} \frac{(n+k+k^*)_{m-k^*} (k)_{k^*}}{(n)_m}.
\end{equation}

The expected value of the number of new boxes in the $m$-sample conditioned to the basic sample, which provides the Bayes estimator for $K_m$ under quadratic loss function, results
\begin{equation}
\label{meanKm}
E(K_m |K_n=k)=\frac{(k)_{n+m}}{(n+\gamma)_{m}}\sum_{k^* =0}^{m} {m \choose k^*} \frac{k^*}{(k+k^*)_n}\frac{(k-\gamma)_{k^*}(\gamma +n -k)_{m-k^*}}{(n)_m}.
\end{equation}

Notice that for $m \rightarrow \infty$, (\ref{bioalpha}) agrees with the posterior result for the total number of boxes obtained in Gnedin (2009). In fact, by standard asymptotic $\Gamma(n+a)/\Gamma(n+b) \sim n^{a-b}$ and recalling the definition of rising and falling factorials in terms of Gamma function 
$$(a)_{b\uparrow} =\frac{\Gamma(a+b)}{\Gamma(a)}   \mbox{       and      }  (a)_{b \downarrow}=\frac{\Gamma(a+1)}{\Gamma(a - b +1)},$$
equation (\ref{bioalpha}) reduces to 
\begin{equation}
\label{bioalphaLIM}
Pr(K=k^* | K_n=k)=\frac{(n-1)!}{(k-1)!} (\gamma +n-k)_{k\uparrow} \frac{(k- \gamma)_{k^* \uparrow} \Gamma(k+k^*)}{\Gamma(k^* +1)\Gamma(n+k+k^*)}.
\end{equation}

The posterior distribution obtained in Gnedin (2009), in terms of  $\varkappa=k+k^*$, is
$$
Pr(\Xi=\varkappa|K_n =k)=\frac{(n-1)!}{(k-1)!(\varkappa + n -1)!} \prod_{i=1}^{k-1}(\varkappa - i) \prod_{j=1}^{k}(\gamma + n -j)\prod_{l=k}^{\varkappa - 1} (l -\gamma)
$$
for $1 \leq k \leq n$, $\varkappa \geq k$ and may be re-written as
$$
=\frac{(n-1)!}{(k-1)!(\varkappa + n - 1)!} (\varkappa -1)_{k-1 \downarrow} (k-\gamma)_{\varkappa -k \uparrow} (\gamma + n -k)_{k \uparrow}.
$$

By substitution $\varkappa=k+k^*$ the result easily follows 
$$
Pr(K=k^* | K_n=k)=\frac{(n-1)!}{(k-1)!} (\gamma +n -k)_{k \uparrow}\frac{\Gamma(k^*+k)}{\Gamma(k+k^*+n)}\frac{(k-\gamma)_{k^*}}{k^*!}.
$$

The corresponding expected value results
$$
E(K| K_n=k)=\frac{(n-1)!}{(k-1!)} (n+\gamma-1)_{k \downarrow} \sum_{k^*=0^\infty} \frac{1}{(k^*-1)!} \frac{(k -\gamma)_{k^*}}{(k+k^*)_n}
$$
which expressed in terms of $\varkappa=k+k^*$ yields
$$
E(\Xi| K_n=k)= \frac{(n-1)!}{(k-1)!}\frac{(n + \gamma-1)_{k \downarrow}}{(k -\gamma-1)!}\sum_{\varkappa=k}^{\infty} \frac{(\varkappa-\gamma -1)!}{(\varkappa - k -1)!}\frac{1}{(\varkappa)_n}.
$$ 

The distribution of the number $S$ of balls in the new $m$-sample which belong to new boxes, given the number of boxes in the basic sample $K_n$ and the number of new boxes $K_m$, follows by an application of Eq. 12 in Lijoi et al. (2008)  
\begin{equation}
\label{eq12aoap}
Pr(S=s|K_m=k^*, K_n=k)={{s-1 \choose k^*-1}{n+k+m-s-1 \choose m-s}} \slash {{n+k+m-1 \choose m-k^*}}
\end{equation}

By Proposition 2. in Lijoi et al. (2008) the mean number of balls  in the subsequent $m$ sample in given by
\begin{equation}
\label{meanS}
E(S|K_n=k)=m\frac{V_{n+1,k+1}}{V_{n,k}}=m\frac{k}{n}\frac{(k-\gamma)}{(n+\gamma)}
\end{equation}
and by Proposition 4 in Lijoi et al. (2008) (cfr. also Corollary 10, in Cerquetti, 2008) the probability that the $m$ new balls don't occupy a subset of $(k-r)$ old boxes arises from (\ref{oldandnew}) by summing over the ways to choose $s$ balls from the $m$ of the new group, by summing over the ways to partition $s$ balls in a subset of $k^*$ boxes, and over the ways to allocate $m-s$ balls in at most $r$ old boxes and is equal to 
\begin{equation}
\label{prop4}
\sum_{k^*=1}^m \frac{(k)_k^*}{(n)_m}\frac{(k -\gamma)_k^*(\gamma +n -k)_{m-k^*}}{(\gamma+n)_m}{m \choose k^*} \frac{1}{(r+\sum_j n_j +m)_{k^* -m}}.
\end{equation}

The conditional Gibbs structure characterizing Gnedin's model as from Proposition 3. in Lijoi et al. (2008), which may be obtained by the operation of {\it deletion of the first $k$ classes} (Pitman, 2003) as clarified in Cerquetti (2008, Prop. 12) will be as follows
$$
p(s_1,\dots, s_{k^*}|{\bf m},{\bf n})
= \frac{\Gamma(k+k^*) \Gamma(k+k^*-\gamma) \Gamma(\gamma + n +m-k -k^*)}{\sum_{k^*=1}^m  {s \choose k^*} \Gamma(k^*)_k {\Gamma(s)} \Gamma(k+k^*-\gamma) \Gamma(\gamma + n +m-k -k^*) } \prod_{j=1}^{k^*} s_j!.
$$

Finally a Bayesian nonparametric estimator for the probability of ball $n+m+1$th to fall in a new box given $K_n=k$ is readily derived by eq. (6) in Lijoi et al. (2007)
$$
\hat{D}^{n,k}_m=\sum_{k^*=0}^m \frac{(k)_{k^* +1}}{(n)_{m+1}} \frac{(k -\gamma)_{k^* +1 } (\gamma + n -k)_{m-k^*}}{(\gamma +n)_{m+1}} {m \choose k^*} \frac{m+n+k^*-1!}{(n-1)!}.
$$

\section{Appendix}

For $n=0,1,2,\dots,$ and arbitrary real $x$ and $h$, let $(x)_{n\uparrow h}$ denote the $n$th factorial power of $x$ with increment $h$ (also called generalized {\it rising} factorial)
\begin{equation}
\label{factdef}
(x)_{n \uparrow h}:= x(x+h)\cdots(x+(n-1)h)=\prod_{i=0}^{n-1}(x+ih)=h^n(x/h)_{n \uparrow},
\end{equation}
where $(x)_{n \uparrow}$ stands for $(x)_{n\uparrow 1}$, $(x)_{h \uparrow 0}=x^h$ and $(x)_{0 \uparrow h}=1$, and for which the following multiplicative law holds 
\begin{equation}
\label{multiplicative}
(x)_{n+r \uparrow h}=(x)_{n\uparrow h} (x +n h)_{r \uparrow h}.
\end{equation}

From e.g. Normand (2004, cfr. eq. 2.41 and 2.45) a binomial formula also holds, namely
\begin{equation}
\label{bino}
(x+y)_{n \uparrow h}=\sum_{k=0}^n {n \choose k} (x)_{k \uparrow h} (y)_{n-k \uparrow h},
\end{equation}
as well as a generalized version of the multinomial theorem, i.e.
\begin{equation}
\label{multi}
(\sum_{j=1}^p z_j)_{n \uparrow h}= \sum_{n_j \geq 0, \sum n_j=n} \frac{n!}{n_1!\cdots n_p!} \prod_{j=1}^p (z_j)_{n_j \uparrow h}.
\end{equation}

We recall the notion of {\it generalized Stirling numbers}, (for a comprehensive treatment see Hsu and Shiue, 1998; see also Pitman, 2006). For arbitrary distinct reals $\eta$ and $\beta$, these are the connection coefficients $S_{n,k}^{\eta, \beta}$ defined by
$$
(x)_{n \downarrow \eta}= \sum_{k=0}^n S_{n,k}^{\eta, \beta} (x)_{k \downarrow \beta} 
$$
where $(x)_{n \downarrow h}$ are generalized {\it falling} factorials and $(x)_{n \downarrow -h}=(x)_{n \uparrow h}$. Hence for $\eta=-1$, $\beta=-\alpha$, and $\alpha \in (-\infty, 1)$, $S_{n,k} ^{-1, -\alpha}$ is defined by
\begin{equation}
\label{unoalpha}
(x)_{n \uparrow 1}=\sum_{k=0}^{n} S_{n,k} ^{-1, -\alpha} (x)_{k \uparrow \alpha},
\end{equation}
or specializing partial Bell polynomials as follows
\begin{equation}
\label{bellalpha}
B_{n,k}((1-\alpha)_{\bullet-1 \uparrow})=\sum_{\{A_1,\dots, A_k\}\in \mathcal{P}_{[n]}^k }\prod_{i=1}^k (1-\alpha)_{n_i-1\uparrow}=\frac{n!}{k!}\sum_{(n_1,\dots,n_k)}\prod_{i=1}^k \frac{(1-\alpha)_{n_i-1\uparrow}}{n_i!}=S_{n,k}^{-1,-\alpha}.
\end{equation}

Referring to formulas in Lijoi et al. (2007, 2008) it is convienent to recall that their treatment is in terms of {\it generalized factorial coefficients}, which are the connection coefficients  $\mathcal{C}^\alpha_{n,k}$ defined by
$
(\alpha y)_{n\uparrow 1}=\sum_{k=0}^{n} \mathcal{C}^\alpha_{n,k}(y)_{k\uparrow 1},
$
(cfr. Charalambides, 2005). 
From (\ref{factdef}) and (\ref{unoalpha}), 
if $x=y \alpha$ then
$$
(y \alpha)_{n \uparrow 1}= \sum_{k=0}^{n} S_{n, k}^{-1, -\alpha}(y \alpha)_{k \uparrow \alpha}=\sum_{k=0}^n S_{n,k}^{-1, -\alpha} \alpha^k
(y)_{k \uparrow 1},
$$
hence $S_{n,k}^{-1, -\alpha}=\alpha^{-k}{\mathcal{C}_{n,k}^\alpha}$.\\

It is also worth to clarify the relationship between {\it non central} generalized Stirling numbers of the first kind  as defined in Hsu and Shiue (1998) and {\it non central} generalized factorial coefficients as in Charalambides (2005).\\\\
First recall that non central generalized Stirling numbers of the first kind are connection coefficients defined by
$$
(x)_{n\uparrow} = \sum_{k=0}^{n} S_{n,k}^{-1,-\alpha, \gamma} (x -\gamma)_{k \uparrow \alpha}
$$
or by the following convolution relation 
\begin{equation}
S_{n,k}^{-1, -\alpha, \gamma}=\sum_{s=0}^{n} {n \choose k} S_{s,k}^{-1, -\alpha} S_{n-s,0} ^{-1, -\alpha, \gamma}.
\end{equation}

Since as a convention we assume
$
S_{s,k}^{-1, -\alpha}=0  \mbox{ for  $s <k$ } 
$
and it is known that 
$
S_{n-s,0}^{-1, -\alpha, \gamma}=(\gamma)_{n-s \uparrow}
$
then
\begin{equation}
\label{convo}
S_{n,k}^{-1, -\alpha, \gamma}= \sum_{s=k}^{n} {n \choose s} S_{s,k}^{-1, -\alpha} (-\gamma)_{n-s \uparrow 1},
\end{equation}
Now, for parameter of non centrality $-\gamma$, and $x=ya$
$$
(y\alpha -\gamma)_{n\uparrow} = \sum_{k=0}^{n} S_{n,k}^{-1,-\alpha, -\gamma} (y\alpha)_{k \uparrow \alpha}
=\sum_{k=0}^n \left[ \sum_{s=k}^n {n \choose s} S_{s,k}^{-1, -\alpha} (-\gamma)_{n-s \uparrow}\right] (y\alpha )_{k \uparrow \alpha}.
$$

Then exploiting the relation between central generalized Stirling numbers and central generalized factorial coefficients
$$
(y \alpha -\gamma)_{n \uparrow}=\sum_{k=0}^n \left[ \alpha^{-k}\sum_{s=k}^n {n \choose s} {C_{s,k}^{\alpha}} (-\gamma)_{n-s \uparrow}\right] (y\alpha)_{k \uparrow \alpha}
$$
and the definition of non central factorial coefficients as in Charalambides (2005)  follows
$$
(y\alpha -\gamma)_{n \uparrow} =\sum_{k=0}^{n}\alpha^{-k} C_{s,k}^{\alpha, \gamma} (ya)_{k \uparrow \alpha}=
\sum_{k=0}^{n} C_{s,k}^{\alpha, \gamma} (y)_{k \uparrow }.
$$

\section*{References}
\newcommand{\bibu}{\item \hskip-1.0cm}
\begin{list}{\ }{\setlength\leftmargin{1.0cm}}







\bibu \textsc {Cerquetti, A.} (2008) Generalized Chinese restaurant construction of exchangeable Gibbs partitions and related results. arXiv:0805.3853v1 [math.PR]

\bibu  \textsc {Charalambides, C. A.} (2005) {\it Combinatorial Methods in Discrete Distributions}. Wiley, Hoboken NJ.





\bibu \textsc{Gnedin, A.} (2009) A species sampling model with finitely many types. arXiv:0910.1988v1 [math.PR]

\bibu \textsc{Gnedin, A. and Pitman, J. } (2006) {Exchangeable Gibbs partitions  and Stirling triangles.} {\it Journal of Mathematical Sciences}, 138, 3, 5674--5685. 



\bibu \textsc{Hsu, L. C, \& Shiue, P. J.} (1998) A unified approach to generalized Stirling numbers. {\it Adv. Appl. Math.}, 20, 366-384.

\bibu \textsc{Lijoi, A., Mena, R. and Pr\"unster, I.} (2007) Bayesian nonparametric estimation of the probability of discovering new species  {\it Biometrika}, 94, 769--786.

\bibu \textsc{Lijoi, A., Pr\"unster, I. and Walker, S.G.} (2008) Bayesian nonparametric estimator derived from conditional Gibbs structures. {\it Annals of Applied Probability}, 18, 1519--1547.

\bibu \textsc {Normand, J.M.} (2004) Calculation of some determinants using the $s$-shifted factorial. {\it J. Phys. A: Math. Gen.} 37, 5737-5762.


\bibu \textsc{Pitman, J.} (1996) Some developments of the Blackwell-MacQueen urn scheme. In T.S. Ferguson, Shapley L.S., and MacQueen J.B., editors, {\it Statistics, Probability and Game Theory}, volume 30 of {\it IMS Lecture Notes-Monograph Series}, pages 245--267. Institute of Mathematical Statistics, Hayward, CA.

\bibu \textsc{Pitman, J.} (2003) {Poisson-Kingman partitions}. In D.R. Goldstein, editor, {\it Science and Statistics: A Festschrift for Terry Speed}, volume 40 of Lecture Notes-Monograph Series, pages 1--34. Institute of Mathematical Statistics, Hayward, California.

\bibu \textsc{Pitman, J.} (2006) {\it Combinatorial Stochastic Processes}. Ecole d'Et\'e de Probabilit\'e de Saint-Flour XXXII - 2002. Lecture Notes in Mathematics N. 1875, Springer.

\bibu \textsc{Pitman, J. and Yor, M.} (1997) The two-parameter Poisson-Dirichlet distribution derived from a stable subordinator. {\it Ann. Probab.}, 25:855--900.







\end{list}
\end{document}